\documentclass[11pt]{amsart}

\begin{document}

\begin{center}
{\large \bf THE UNIVERSAL MUMFORD CURVE, 
AND ITS ABELIAN DIFFERENTIALS AND PERIODS 
IN ARITHMETIC FORMAL GEOMETRY}
\footnote{Research supported by the JSPS Grant-in-Aid for 
Scientific Research No. 20K03516.} 
\end{center}
\vspace{2ex}

\begin{center}
{\large \sc By Takashi Ichikawa} 
\end{center}
\vspace{2ex}

\hspace{5cm} \hrulefill \hspace{5cm}
\vspace{3ex}



\noindent
\begin{small} 
{\it Abstract.} 
We construct the universal Mumford curve of given genus as a family of Mumford curves 
over the deformation space of degenerate curves 
in the category of arithmetic formal geometry. 
Furthermore, we give explicit formulas of abelian differentials and their periods of 
the universal Mumford curve.   
\end{small}
\vspace{2ex}

{\bf 1. Introduction.} 
This paper is a revision of \cite{I8} with some corrections and extensions. 
As an analog of the Schottky uniformization theory of Riemann surfaces \cite{S}, 
Mumford \cite{Mu1} established a uniformization theory of Mumford curves 
defined as stable curves over complete local rings whose special fibers are degenerate, 
i.e., consisting of (may be self-intersecting) projective lines. 
There were many researches on this theory, 
and in recent works of Ulirsch \cite{U} and Poineau-Turchetti \cite{PT}, 
based on tropical geometry and Berkovich geometry over ${\mathbb Z}$ \cite{P1, P2}, 
a {\it universal Mumford curve} was given unifying Riemann surfaces and 
Mumford curves over $p$-adic fields. 
Considering the Raynaud correspondence \cite{R} 
between formal geometry and rigid analytic geometry, 
one can consider the problem of constructing such a curve 
in the category of formal geometry over ${\mathbb Z}$. 

Our purpose of this paper is to give a solution of the problem 
in terms of {\it arithmetic formal geometry} over ${\mathbb Z}$ 
as constructing the universal Mumford curve which gives rise to all Mumford curves, 
and to all Riemann surfaces close to degenerate complex curves. 
The $n$-marked universal Mumford curve of genus $g$ is defined over the deformation space  ${\mathcal A}_{g, n}$ of all $n$-marked degenerate curves of genus $g$ 
whose analytic counterpart is the quotient of the extended Schottky space  $\overline{S}_{g,n}$ (cf. \cite{GH}) by the outer automorphism group ${\rm Out}(F_{g})$ 
of the free group $F_{g}$ of rank $g$. 
Comparing with the result of \cite{PT}, 
our theory is directly connected with the deformation theory of degenerate curves, 
and gives explicit formulas of abelian differentials and period maps 
for degenerating curves as is stated below. 
Furthermore, the evaluation of the universal Mumford curve provides tools of studying  arithmetic properties and asymptotic behaviors of functions defined for families of curves, 
e.g., Teichm\"{u}ller modular forms and their variants (cf. \cite{I1, I4, I5, I6}). 

We construct the universal Mumford curve using generalized Tate curves defined in \cite{I1} 
which are universal families of curves with given degeneration data 
described by the moduli and deformation parameters over ${\mathbb Z}$. 
In this paper, based on results of \cite{I2, I3} for Teichm\"{u}ller groupoids with application to Teichm\"{u}ller's lego game of Grothendieck \cite{Gr} (see also \cite{BaK1, BaK2, MoS, NS}), 
we obtain the universal Mumford curve by gluing generalized Tate curves 
comparing the associated parameters. 

Another purpose of this paper is to give an application of the universal Mumford curve 
by describing explicit formulas of its abelian differentials and their period integrals. 
Our formulas unify results of Schottky \cite{S}, Bainbridge-M\"{o}ller \cite{BM}, 
Hu-Norton \cite{HN}, Manin-Drinfeld \cite{MD}, de Shalit \cite{dS1, dS2} and others 
for families of Riemann surfaces and of Mumford curves 
in terms of arithmetic formal geometry (cf. \cite{I7}). 
Furthermore, we construct the (generalized) Jacobian of the universal Mumford curve 
by gluing Jacobians of generalized Tate curves 
based on results of Mumford \cite{Mu2} and Faltings-Chai \cite{FC}. 
Alexeev \cite{A}, Caporaso \cite{C}, Esteves \cite{E}, Oda-Seshadri \cite{OS} and 
Simpson \cite{Si} gave moduli theoretic constructions of universal compactified Jacobians 
of relative dimension $g$ which are seen to be associated with 
special polyhedral cone decompositions of ${\mathbb R}^{g}$ (cf. \cite{ChPS}). 
For general ${\mathbb Z}^{g}$-admissible polyhedral cone decompositions of ${\mathbb R}^{g}$, 
our result will imply explicit construction of compactified Jacobians which are 
{\it semi-global} as they are defined on ${\mathcal A}_{g}$ and on 
$\overline{S}_{g} / {\rm Out}(F_{g})$ by the analytic extension. 
\vspace{4ex}

{\bf 2. Generalized Tate curve.}  
\vspace{2ex}

{\bf 2.1. Schottky uniformization.} \ 
A Schottky group $\Gamma$ of rank $g$ is defined as a free group with generators 
$\gamma_{i} \in PGL_{2}({\mathbb C})$ $(i = 1,..., g)$ which map Jordan curves  
$C_{i} \subset {\mathbb P}^{1}_{\mathbb C} = {\mathbb C} \cup \{ \infty \}$ 
to other Jordan curves $C_{-i} \subset {\mathbb P}^{1}_{\mathbb C}$ 
with orientation reversed, 
where $C_{\pm 1},..., C_{\pm g}$ with their interiors are mutually disjoint. 
Each element $\gamma \in \Gamma - \{ 1 \}$ is conjugate 
to an element of $PGL_{2}({\mathbb C})$ sending $z$ to $\beta_{\gamma} z$ 
for some $\beta_{\gamma} \in {\mathbb C}^{\times}$ with $|\beta_{\gamma}| < 1$ 
which is called the {\it multiplier} of $\gamma$. 
Therefore, one has 
$$
\frac{\gamma(z) - \alpha_{\gamma}}{z - \alpha_{\gamma}} = 
\beta_{\gamma} \frac{\gamma(z) - \alpha'_{\gamma}}{z - \alpha'_{\gamma}} 
$$
for some element $\alpha_{\gamma}, \alpha'_{\gamma}$ of 
${\mathbb P}^{1}_{\mathbb C}$ 
called the {\it attractive}, {\it repulsive} fixed points of $\gamma$ 
respectively. 
Then the discontinuity set 
$\Omega_{\Gamma} \subset {\mathbb P}^{1}_{\mathbb C}$ 
under the action of $\Gamma$ has a fundamental domain $D_{\Gamma}$ 
which is given by the complement of the union of the interiors of $C_{\pm 1},..., C_{\pm g}$. 
The quotient space $R_{\Gamma} = \Omega_{\Gamma}/\Gamma$ 
is a (compact) Riemann surface of genus $g$ 
which is called {\it Schottky uniformized} by $\Gamma$ (cf. \cite{S}). 
Furthermore, by a result of Koebe, 
every Riemann surface of genus $g$ can be represented in this manner. 
\vspace{2ex}

{\bf 2.2. Generalized Tate curve.} \ 
A (marked) curve is called {\it degenerate} if it is a stable (marked) curve and 
the normalization of its irreducible components are all projective (marked) lines. 
Then the dual graph $\Delta = (V, E, T)$ of a stable marked curve is a collection  
of 3 finite sets $V$ of vertices, $E$ of edges, $T$ of tails 
and 2 boundary maps 
$$
b : T \rightarrow V, 
\ \ b : E \longrightarrow \left( V \cup \{ \mbox{unordered pairs of elements of $V$} \} \right) 
$$
such that the geometric realization of $\Delta$ is connected 
and that $\Delta$ is {\it stable}, namely its each vertex has at least $3$ branches. 
The number of elements of a finite set $X$ is denoted by $\sharp X$, 
and a (connected) stable graph $\Delta = (V, E, T)$ is called {\it of $(g, n)$-type} 
if ${\rm rank}_{\mathbb Z} H_{1}(\Delta, {\mathbb Z}) = g$, $\sharp T = n$. 
Then under fixing a bijection $\nu : T \stackrel{\sim}{\rightarrow} \{ 1, ... , n \}$, 
which we call a numbering of $T$, 
$\Delta = (V, E, T)$ becomes the dual graph of a degenerate $n$-marked curve of genus $g$ 
such that each tail $h \in T$ corresponds to the $\nu(h)$th marked point. 
In particular, a stable graph without tail is the dual graph of 
a degenerate (unmarked) curve by this correspondence. 
If $\Delta$ is trivalent, i.e. any vertex of $\Delta$ has just $3$ branches, 
then a degenerate $\sharp T$-marked curve with dual graph $\Delta$ 
is maximally degenerate. 
An {\it orientation} of a stable graph $\Delta = (V, E, T)$ means 
giving an orientation of each $e \in E$. 
Under an orientation of $\Delta$, 
denote by $\pm E = \{ e, -e \ | \ e \in E \}$ the set of oriented edges, 
and by $v_{h}$ the terminal vertex of $h \in \pm E$ (resp. the boundary vertex of $h \in T)$. 
For each $h \in \pm E$, 
denote by $|h| \in E$ the edge $h$ without orientation. 

Let $\Delta = (V, E, T)$ be a stable graph. 
Fix an orientation of $\Delta$, 
and take a subset ${\mathcal E}$ of $\pm E \cup T$ 
whose complement ${\mathcal E}_{\infty}$ satisfies the condition that 
$$
\pm E \cap {\mathcal E}_{\infty} \cap \{ -h \ | \ h \in {\mathcal E}_{\infty} \} 
\ = \ 
\emptyset, 
$$ 
and that $v_{h} \neq v_{h'}$ for any distinct $h, h' \in {\mathcal E}_{\infty}$. 
We attach variables $x_{h}$ for $h \in {\mathcal E}$ and $y_{e} = y_{-e}$ for $e \in E$ 
which we call {\it moduli parameters} and {\it deformation parameters} respectively. 
Let $R_{\Delta}$ be the ${\mathbb Z}$-algebra generated by $x_{h}$ $(h \in {\mathcal E})$, 
$1/(x_{e} - x_{-e})$ $(e, -e \in {\mathcal E})$ and $1/(x_{h} - x_{h'})$ 
$(h, h' \in {\mathcal E}$ with $h \neq h'$ and $v_{h} = v_{h'})$, 
and let 
$$ 
A_{\Delta} \ = \ R_{\Delta} [[y_{e} \ (e \in E)]], \ \ 
B_{\Delta} \ = \ A_{\Delta} \left[ \prod_{e \in E} y_{e}^{-1} \right]. 
$$ 

According to \cite[Section 2]{I1}, 
we construct the {\it universal Schottky group} $\Gamma_{\Delta}$ 
associated with oriented $\Delta$ and ${\mathcal E}$ as follows. 
For $h \in \pm E$, put 
\begin{eqnarray*}
\phi_{h} 
& = & 
\left( \begin{array}{cc} x_{h} & x_{-h} \\ 1 & 1 \end{array} \right) 
\left( \begin{array}{cc} 1 & 0 \\ 0 & y_{h} \end{array} \right) 
\left( \begin{array}{cc} x_{h} & x_{-h} \\ 1 & 1 \end{array} \right)^{-1} 
\\ 
& = & 
\frac{1}{x_{h} - x_{-h}} 
\left\{ \left( \begin{array}{cc} x_{h} & - x_{h} x_{-h} \\ 1 & - x_{-h} \end{array} \right) - 
\left( \begin{array}{cc} x_{-h} & - x_{h} x_{-h} \\ 1 & - x_{h} \end{array} \right) y_{h} \right\}, 
\end{eqnarray*} 
where $x_{h}$ (resp. $x_{-h})$ means $\infty$ 
if $h$ (resp. $-h)$ belongs to ${\mathcal E}_{\infty}$, 
which gives an element of $PGL_{2}(B_{\Delta}) = GL_{2}(B_{\Delta})/B_{\Delta}^{\times}$ 
denoted by the same symbol. 
Its attractive (resp. repulsive) fixed points are $x_{h}$ (resp. $x_{-h}$), 
and its multiplier is $y_{h}$, namely 
$$
\frac{\phi_{h}(z) - x_{h}}{z - x_{h}} 
\ = \ 
y_{h} \frac{\phi_{h}(z) - x_{-h}}{z - x_{-h}} 
\ \ \left( z \in {\mathbb P}^{1} \right). 
$$
For any reduced path $\rho = h(1) \cdot h(2) \cdots h(l)$ 
which is the product of oriented edges $h(1), ... ,h(l)$ 
such that $v_{h(i)} = v_{-h(i+1)}$ and $h(i) \neq - h(i+1)$, 
one can associate an element $\rho^{*}$ of $PGL_{2}(B_{\Delta})$ 
having reduced expression $\phi_{h(l)} \phi_{h(l-1)} \cdots \phi_{h(1)}$. 
Fix a base vertex $v_{b}$ of $V$, 
and consider the fundamental group 
$\pi_{1} (\Delta, v_{b})$ which is a free group 
of rank $g = {\rm rank}_{\mathbb Z} H_{1}(\Delta, {\mathbb Z})$. 
Then the correspondence $\rho \mapsto \rho^{*}$ 
gives an injective anti-homomorphism 
$\pi_{1} (\Delta, v_{b}) \rightarrow PGL_{2}(B_{\Delta})$ 
whose image is denoted by $\Gamma_{\Delta}$. 

It is shown in \cite[Section 3]{I1} and \cite[1.4]{I2} 
(see also \cite[Section 2]{IhN} when $\Delta$ is trivalent and has no loop) 
that for any stable graph $\Delta$ of $(g, n)$-type, 
there exists a stable $n$-marked curve $C_{\Delta}$ of genus $g$ over $A_{\Delta}$ 
which satisfy the following properties: 

\begin{itemize}

\item[(P1)] 
The closed fiber $C_{\Delta} \otimes_{A_{\Delta}} R_{\Delta}$ of $C_{\Delta}$ 
obtained by substituting $y_{e} = 0$ $(e \in E)$ 
becomes the degenerate marked curve over $R_{\Delta}$ with dual graph $\Delta$ which is 
obtained from the collection of $P_{v} := {\mathbb P}^{1}_{R_{\Delta}}$ $(v \in V)$ 
by identifying the points $x_{e} \in P_{v_{e}}$ and $x_{-e} \in P_{v_{-e}}$ ($e \in E$), 
where $x_{h}$ denotes $\infty$ if $h \in {\mathcal E}_{\infty}$. 

\item[(P2)] 
$C_{\Delta}$ gives rise to a universal deformation 
of degenerate marked curves with dual graph $\Delta$. 
More precisely, $C_{\Delta}$ satisfies the following: 
Let $R$ be a noetherian and normal complete local ring with residue field $k$, 
and let $C$ be a marked Mumford curve over $R$, 
namely a stable marked curve over $R$ with nonsingular generic fiber 
such that the closed fiber $C \otimes_{R} k$ is a degenerate marked curve 
in which all double points and marked points are $k$-rational. 
Then under the condition that $\Delta$ is the dual graph of $C \otimes_{R} k$, 
there exists a ring homomorphism $A_{\Delta} \rightarrow R$ 
giving $C_{\Delta} \otimes_{A_{\Delta}} R \cong C$.  

\item[(P3)] 
$C_{\Delta} \otimes_{A_{\Delta}} B_{\Delta}$ is smooth over $B_{\Delta}$ 
and is Mumford uniformized (cf. \cite{Mu1}) by $\Gamma_{\Delta}$. 

\item[(P4)] 
Take $x_{h}$ $(h \in {\mathcal E})$ as complex numbers such that $x_{e} \neq x_{-e}$ 
and that $x_{h} \neq x_{h'}$ if $h \neq h'$ and $v_{h} = v_{h'}$, 
and take $y_{e}$ $(e \in E)$ as sufficiently small nonzero complex numbers. 
Then $C_{\Delta}$ becomes a marked Riemann surface which is 
Schottky uniformized by the Schottky group $\Gamma$ over ${\mathbb C}$ 
obtained from $\Gamma_{\Delta}$. 

\end{itemize}
We call  $C_{\Delta}$ the {\it generalized Tate curve} associated with $\Delta$, 
and review its construction given in \cite[Theorem 3.5]{I1}. 
Let $T_{\Delta}$ be the tree obtained as the universal cover of $\Delta$, 
and denote by ${\mathcal P}_{T_{\Delta}}$ be the formal scheme 
as the union of copies of ${\mathbb P}^{1}_{A_{\Delta}}$ indexed by vertices of $T_{\Delta}$ 
under the $B_{\Delta}$-isomorphism by $\phi_{e}$ $(e \in E)$. 
Then it is shown in \cite[Theorem 3.5]{I1} that 
$C_{\Delta}$ is the formal scheme theoretic quotient of ${\mathcal P}_{T_{\Delta}}$ by $\Gamma_{\Delta}$. 
\vspace{4ex}

{\bf 3. Universal Mumford curve.}
\vspace{2ex}

{\bf 3.1. Comparison of deformations.} \ 
Let $\Delta = (V, E, T)$ be a stable graph which is not trivalent. 
Then there exists a vertex $v_{0} \in V$ which has at least $4$ branches. 
Take two elements $h_{1}, h_{2}$ of ${\mathcal E}$ such that 
$h_{1} \neq h_{2}$ and $v_{h_{1}} = v_{h_{2}} = v_{0}$, 
and let $\Delta' = (V', E', T')$ be a stable graph obtained from $\Delta$ 
by replacing $v_{0}$ with an oriented (nonloop) edge $h_{0}$ such that 
$v_{h_{1}} = v_{h_{2}} = v_{h_{0}}$ and that $v_{h} = v_{-h_{0}}$ 
for any $h \in \pm E \cup T - \{ h_{1}, h_{2} \}$ with $v_{h} = v_{0}$. 
Put $e_{i} = |h_{i}|$ for $i = 0, 1, 2$. 
Then we have the following identifications: 
$$
V = V' - \{ v_{-h_{0}} \} \ (\mbox{in which $v_{0} = v_{h_{0}}$}), \ E = E' - \{ e_{0} \}, \ T = T'. 
$$

{\sc Theorem 3.1.} 
\begin{it} 

{\rm (1)} 
The generalized Tate curves $C_{\Delta}$ and $C_{\Delta'}$ 
associated with $\Delta$ and $\Delta'$ respectively are isomorphic 
over $R_{\Delta'} [[ s_{e} \ (e \in E') ]] [s_{e_{0}}^{-1}]$, 
where 
$$
\frac{x_{h_{1}} - x_{h_{2}}}{s_{e_{0}}}, \ \ 
\frac{y_{e_{i}}}{s_{e_{0}} s_{e_{i}}} \ (\mbox{$i = 1, 2$ with $h_{i} \not\in T$}), \ \
\frac{y_{e}}{s_{e}} \ (e \in E - \{ e_{1}, e_{2} \}) 
$$
belong to $(A_{\Delta'})^{\times}$ if $h_{1} \neq - h_{2}$, 
and 
$$
\frac{x_{h_{1}} - x_{h_{2}}}{s_{e_{0}}}, \ \ 
\frac{y_{e}}{s_{e}} \ (e \in E) 
$$
belong to $(A_{\Delta'})^{\times}$ if $h_{1} = - h_{2}$. 

{\rm (2)} 
The assertion (1) holds in the category of complex geometry 
when $x_{h_{1}} - x_{h_{2}}, y_{e}$ and $s_{e}$ are taken to be 
sufficiently small complex numbers with $x_{h_{1}} \neq x_{h_{2}}, s_{e_{0}} \neq 0$. 

\end{it}
\vspace{2ex}

{\it Proof.} 
First, we prove the assertion (1). 
Denote by $t_{h}$ the moduli parameters of $C_{\Delta'}$ corresponding to 
$h \in \pm E' \cup T'$.  
Then by \cite[Lemma 1.2]{I1}, 
$\phi_{-h_{0}}(t_{h_{1}}) - \phi_{-h_{0}}(t_{h_{2}})$ belongs to 
$s_{e_{0}} \cdot (A_{\Delta'})^{\times}$, 
and hence 
$$
C_{\Delta'} \otimes_{A_{\Delta'}} 
\left( R_{\Delta'} [[ s_{e} ]] [s_{e_{0}}^{-1}] \right)
$$ 
gives a universal deformation of a universal degenerate curve with dual graph $\Delta$. 
Then by the universality of generalized Tate curves, 
there exists an injective homomorphism 
$A_{\Delta} \hookrightarrow R_{\Delta'} [[ s_{e} ]] [s_{e_{0}}^{-1}]$
which gives rise to an isomorphism $C_{\Delta} \cong C_{\Delta'}$. 
Under this homomorphism, 
\begin{eqnarray*}
\lefteqn{
\left( P_{v_{-h_{0}}}; \phi_{-h_{0}}(t_{h_{1}}), \phi_{-h_{0}}(t_{h_{2}}), 
t_{h} \ (v_{h} = v_{-h_{0}}, h \neq h_{0}) \right) 
} 
\\ 
& \cong & 
\left( P_{v_{0}}; x_{h_{1}}, x_{h_{2}}, x_{h} \ (v_{h} = v_{0}, h \neq h_{1}, h_{2}) \right), 
\end{eqnarray*} 
and hence $x_{h_{1}} - x_{h_{2}} \in s_{e_{0}} \cdot (A_{\Delta'})^{\times}$. 
Furthermore, when $h_{1} \neq -h_{2}$, 
the deformation parameters of 
$$
\left( P_{v_{-h_{0}}}; \phi_{-h_{0}}(t_{h_{1}}), \phi_{-h_{0}}(t_{h_{2}}), 
t_{h} \ (v_{h} = v_{-h_{0}}, h \neq h_{0}) \right) 
$$
corresponding to $h_{i} \cdot (-h_{0})$ $(i = 1, 2)$ are $y_{h_{i}}$, 
and hence by \cite[Proposition 1.3]{I1}, 
$y_{h_{i}} \in \left( s_{h_{0}} \cdot s_{h_{i}} \right) \cdot  (A_{\Delta'})^{\times}$. 
When $h_{1} = - h_{2}$, 
the deformation parameters of 
$$
\left( P_{v_{-h_{0}}}; \phi_{-h_{0}}(t_{h_{1}}), \phi_{-h_{0}}(t_{h_{2}}), 
t_{h} \ (v_{h} = v_{-h_{0}}, h \neq h_{0}) \right) 
$$
corresponding to $h_{0} \cdot h_{1} \cdot (-h_{0})$ is $y_{h_{1}}$, 
and hence $y_{h_{1}} \in s_{h_{1}} \cdot  (A_{\Delta'})^{\times}$. 

The assertion (2) follows from the property 2.2 (P4) for generalized Tate curves. 
\ $\square$ 
\vspace{2ex}

{\it Remark} 3.2. 
Let $\Gamma_{\Delta}$ and $\Gamma_{\Delta'}$ denote the universal Schottky groups 
associated with $\Delta$ and $\Delta'$ respectively. 
Since $C_{\Delta} \cong C_{\Delta'}$, 
by \cite[Corollary 4.11]{Mu1}, there exists a conjugation isomorphism 
$\varphi : \Gamma_{\Delta} \stackrel{\sim}{\rightarrow} \Gamma_{\Delta'}$, 
and hence the multiplier of an element of $\Gamma_{\Delta}$ 
and the cross-ratio of the fixed points of four elements of $\Gamma_{\Delta}$ 
are invariant under $\varphi$. 
Therefore, one can calculate the precise formula connecting 
the moduli and deformation parameters for $C_{\Delta}$ and $C_{\Delta'}$. 
For the detail, see the proof of \cite[Theorem 1]{I2}. 
\vspace{2ex}

{\bf 3.2. Construction of the universal Mumford curve.} \ 
For nonnegative integers $g, n$ such that $2g - 2 + n > 0$, 
denote by $\overline{\mathcal M}_{g,n}$ the moduli stack over ${\mathbb Z}$ 
of stable $n$-marked curves of genus $g$ (cf. \cite{DM, K, KM}). 
Then by definition, there exist the universal stable marked curve ${\mathcal C}_{g,n}$ 
over $\overline{\mathcal M}_{g,n}$, 
and the associated curve ${\mathcal C}_{g}$ 
obtained by forgetting marked points on ${\mathcal C}_{g,n}$. 
\vspace{2ex}

{\sc Theorem 3.3.} 
\begin{it} 
There exists a deformation space ${\mathcal A}_{g, n}$ of 
all $n$-marked degenerate curves of genus $g$, 
and an $n$-marked stable curve of genus $g$ over ${\mathcal A}_{g, n}$ 
whose fiber by the canonical morphism ${\rm Spec}(A_{\Delta}) \rightarrow {\mathcal A}_{g, n}$ 
becomes the generalized Tate curve $C_{\Delta}$ 
for each stable graph $\Delta$ of $(g, n)$-type. 
\end{it} 
\vspace{2ex}

{\it Proof.} 
Let $\Delta = (V, E, T)$ be a stable graph of $(g, n)$-type, 
and take a system of coordinates on $P_{v} = {\mathbb P}^{1}_{R_{\Delta}}$ $(v \in V)$ 
such that $x_{h} = \infty$ $(h \in {\mathcal E}_{\infty})$ and that 
$\{ 0, 1 \} \subset P_{v}$ is contained in the set of points 
given by $x_{h}$ $(h \in {\mathcal E}$ with $v_{h} = v)$. 
Under this system of coordinates, 
one has the generalized Tate curve $C_{\Delta}$ whose closed fiber 
$C_{\Delta} \otimes_{A_{\Delta}} R_{\Delta}$ gives a family of degenerate curves 
over the open subspace of 
$$
S_{\Delta} = \left\{ \left. \left( p_{h} \in P_{v_{h}} \right)_{h \in \pm E \cup T} \ \right| \ 
p_{h} \neq p_{h'} \ (h \neq h', v_{h} = v_{h'}) \right\} 
$$ 
defined as $p_{e} \neq p_{-e}$ for nonloop edges $e \in E$. 
Therefore, taking another system of coordinates on $P_{v}$ obtained by mutual changes 
of $0, 1, \infty$ and comparing the associated generalized Tate curves 
with the original $C_{\Delta}$ as in Theorem 3.1, 
$C_{\Delta}$ can be extended over the deformation space of 
all marked degenerate curves with dual graph $\Delta$. 
Since two stable graphs of $(g, n)$-type can be translated by a combination of 
replacements $\Delta \leftrightarrow \Delta'$ given in 3.1, 
one can define a scheme ${\mathcal A}_{g, n}$ obtained by gluing 
${\rm Spec}(A_{\Delta})$ ($\Delta$: stable graphs of $(g, n)$-type) 
along the isomorphism given in Theorem 3.1. 
Then ${\mathcal A}_{g, n}$ is regarded as the deformation space of 
all $n$-marked degenerate curves of genus $g$ over which 
there exists an $n$-marked stable curve of genus $g$ obtained by gluing $C_{\Delta}$. 
\ $\square$  
\vspace{2ex}

{\it Definition} 3.4. 
We call the above $n$-marked stable curve of genus $g$ over ${\mathcal A}_{g, n}$  
the $n$-marked {\it universal Mumford curve} of genus $g$ 
which is the fiber of ${\mathcal C}_{g, n}$ by the canonical morphism 
${\mathcal A}_{g, n} \rightarrow \overline{\mathcal M}_{g, n}$. 
By 2.2 (P2) and (P4), 
one can see that this universal Mumford curve gives rise to 
all $n$-marked Mumford curves of genus $g$, 
and to all $n$-marked Riemann surfaces of genus $g$ close to degenerate curves. 
\vspace{2ex}

{\it Remark} 3.5. 
Gerritzen-Herrlich \cite{GH} introduced the extended Schottky space 
$\overline{S}_{g}$ of genus $g > 1$ 
as the fine moduli space of stable complex curves of genus $g$ with Schottky structure. 
For integers $g, n$ as above, 
one can consider the extended Schottky space $\overline{S}_{g,n}$ 
for stable $n$-marked complex curves of genus $g$ with Schottky structure. 
Then by the result of Koebe referred in 2.1, 
$\overline{S}_{g, n}/{\rm Out}(F_{g})$ becomes a covering space of 
the moduli space of $n$-marked stable complex curves of genus $g$ 
which is also the complex analytic space $\overline{\mathcal M}_{g, n}^{\rm an}$ 
associated with $\overline{\mathcal M}_{g, n}$. 
Furthermore, the $n$-marked universal Mumford curve of genus $g$ can be 
analytically  extended to the universal family of marked stable complex curves 
over $\overline{S}_{g,n}/{\rm Out}(F_{g})$. 
\vspace{4ex}

{\bf 4. Abelian differentials on generalized Tate curves.}
\vspace{2ex}

{\bf 4.1. Abelian differentials.} \ 
Let $\Delta = (V, E, T)$ be a stable graph of $(g, n)$-type, 
where $2g - 2 + n > 0$, and the notation be as in 2.2. 
\vspace{2ex}

{\sc Proposition 4.1.} 
\begin{it}
Let $\phi$ be a product $\phi_{h(1)} \cdots \phi_{h(l)}$ with 
$v_{-h(i)} = v_{h(i+1)}$ $(i = 1,..., l-1)$ 
which is reduced in the sense that $h(i) \neq -h(i+1)$ $(i = 1,..., l-1)$, 
and put $y_{\phi} = y_{h(1)} \cdots y_{h(l)}$. 

{\rm (1)} 
One has 
$\displaystyle \phi(z) - x_{h(1)} \in y_{h(1)} 
\left( R_{\Delta} \left[ z, \prod_{h \in \pm E} (z - x_{h})^{-1} \right] [[ y_{e} \ (e \in E) ]] \right)$. 

{\rm (2)} 
If $a \in A_{\Delta}$ satisfies $a - x_{-h(l)} \in A_{\Delta}^{\times}$, 
then $\phi(a) - x_{h(1)} \in I$. 
Furthermore, if $a' - x_{-h(l)} \in A_{\Delta}^{\times}$, 
then $\phi(a) - \phi(a') \in (a - a') y_{\phi} A_{\Delta}^{\times}$. 

{\rm (3)} 
One has 
$\displaystyle 
\frac{d \phi(z)}{dz} \in y_{\phi} 
\left( R_{\Delta} \left[ \prod_{h \in \pm E} (z - x_{h})^{-1} \right] [[ y_{e} \ (e \in E) ]] \right)$. 
\end{it} 
\vspace{2ex}

{\it Proof.}
Since the assertion (2) is proved in \cite[Lemma 1.2]{I1}, 
we will prove (1) and (3). 
Put 
$$
\phi = \left( \begin{array}{cc} a_{\phi} & b_{\phi} \\ c_{\phi} & d_{\phi} \end{array} \right). 
$$
Since 
$$
\left( \begin{array}{cc} \alpha & - \alpha \beta \\ 1 & - \beta \end{array} \right)
\left( \begin{array}{cc} \gamma & - \gamma \delta \\ 1 & - \delta \end{array} \right) 
= (\gamma - \beta) 
\left( \begin{array}{cc} \alpha & - \alpha \beta \\ 1 & - \delta \end{array} \right), 
$$
$a_{\phi}$, $b_{\phi}$, $c_{\phi}$ and $d_{\phi}$ are elements of $A_{\Delta}$ 
whose constant terms are $x_{h(1)} t$, $- x_{h(1)} x_{-h(l)} t$, $t$ and $- x_{-h(l)} t$ 
respectively,   
where 
$$
t = \frac{\prod_{s = 2}^{l} \left( x_{h(s)} - x_{-h(s-1)} \right)}
{\prod_{s = 1}^{l} \left( x_{h(s)} - x_{-h(s)} \right)} \in A_{\Delta}^{\times}. 
$$
Then $c_{\phi} z + d_{\phi} = t (z - x_{-h(l)}) + \cdots$, 
and hence  
$$
\phi(z) - x_{h(1)} \in 
R_{\Delta} \left[ z, \prod_{h \in \pm E} (z - x_{h})^{-1} \right] [[ y_{e} \ (e \in E) ]]. 
$$
In order to prove (1), we may assume that $l = 1$, 
and then $\phi(z) = x_{h(1)}$ under $y_{h(1)} = 0$. 
Therefore, the assertion (1) holds.  
The assertion (3) follows from 
$$
\frac{d \phi(z)}{dz} = \frac{\det(\phi)}{(c_{\phi} z + d_{\phi})^{2}} 
= \frac{\prod_{s = 1}^{l} \det(\phi_{h(s)})}{(c_{\phi} z + d_{\phi})^{2}} 
= \frac{\prod_{s = 1}^{l} y_{h(s)}}{(c_{\phi} z + d_{\phi})^{2}}, 
$$
and the above calculation. 
\ $\square$ 
\vspace{2ex} 

For a stable graph $\Delta = (V, E, T)$, 
we define abelian differentials on a generalized Tate curve $C_{\Delta}$ 
(cf. \cite{S} and \cite[Corollary of Theorem 3]{MD}). 
Let 
$\Gamma_{\Delta} = {\rm Im} \left( \pi_{1}(\Delta, v_{b}) \rightarrow 
PGL_{2}(B_{\Delta}) \right)$ 
be the universal Schottky group as above. 
Then it is shown in \cite[Proposition 1.3]{I1} that 
each $\gamma \in \Gamma_{\Delta} - \{ 1 \}$ has its attractive (resp. repulsive) 
fixed points $\alpha$ (resp. $\alpha'$) in ${\mathbb P}^{1}_{B_{\Delta}}$ 
and its multiplier $\beta \in \sum_{e \in E} A_{\Delta} \cdot y_{e}$ which satisfy 
$$
\frac{\gamma(z) - \alpha}{z - \alpha} = \beta \frac{\gamma(z) - \alpha'}{z - \alpha'}. 
$$ 
Fix a set $\{ \gamma_{1},..., \gamma_{g} \}$ of generators of $\Gamma_{\Delta}$, 
and for each $\gamma_{i}$, 
denote by $\alpha_{i}$ (resp. $\alpha_{-i}$) its attractive (resp. repulsive) fixed points, 
and by $\beta_{i}$ its multiplier. 
Then under the assumption that there is no element of $\pm E \cap {\mathcal E}_{\infty}$ 
with terminal vertex $v_{b}$, 
for each $i = 1,..., g$, 
we define the associated {\it abelian differential of the first kind} as 
$$
\omega_{i} = 
\sum_{\gamma \in \Gamma_{\Delta} / \left\langle \gamma_{i} \right\rangle} 
\left( \frac{1}{z - \gamma(\alpha_{i})} - \frac{1}{z - \gamma(\alpha_{-i})} \right) dz. 
$$
Assume that 
$$
\{ h \in \pm E \cap {\mathcal E}_{\infty} \ | \ v_{h} = v_{b} \} = \emptyset,  \ \ 
\{ t \in T \ | \ v_{t} = v_{b} \} \neq \emptyset. 
$$ 
Then for each $t \in T$ with $v_{t} = v_{b}$ and $k > 1$, 
we define the associated {\it abelian differential of the second kind} as 
$$
\omega_{t, k} = \left\{ \begin{array}{ll} 
{\displaystyle \sum_{\gamma \in \Gamma_{\Delta}} \frac{d \gamma(z)}{(\gamma(z) - x_{t})^{k}}} 
& (x_{t} \neq \infty), 
\\
\\
{\displaystyle \sum_{\gamma \in \Gamma_{\Delta}} \gamma(z)^{k-2} d \gamma(z)} 
& (x_{t} = \infty). 
\end{array} \right. 
$$
Furthermore, 
put $T_{\infty} = \{ t \in T \cap {\mathcal E}_{\infty} \ | \ v_{t} = v_{b} \}$ 
whose cardinality is $0$ or $1$, 
and take a maximal subtree ${\mathcal T}_{\Delta}$ of $\Delta$, 
and for each $t \in T$, 
take the unique path $\rho_{t} = h(1) \cdots h(l)$ in ${\mathcal T}_{\Delta}$ 
from $v_{t}$ to $v_{b}$, 
and put $\phi_{t} = \phi_{h(l)} \cdots \phi_{h(1)}$. 
Then for each $t_{1}, t_{2} \in T$ with $t_{1} \neq t_{2}$, 
we define the associated {\it abelian differential of the the third kind} as 
$$
\omega_{t_{1}, t_{2}} = \sum_{\gamma \in \Gamma_{\Delta}} 
\left( \frac{d \gamma(z)}{\gamma(z) - \phi_{t_{1}}(x_{t_{1}})} - 
\frac{d \gamma(z)}{\gamma(z) - \phi_{t_{2}}(x_{t_{2}})} \right), 
$$
where $\phi_{t_{i}}(x_{t_{i}}) = \infty$ if $t_{i} \in T_{\infty}$. 
\vspace{2ex}

{\sc Theorem 4.2.} 
\begin{it} 

{\rm (1)} 
For each $i = 1,..., g$, 
$\omega_{i}$ is a regular differential on $C_{\Delta} \otimes_{A_{\Delta}} B_{\Delta}$. 

{\rm (2)} 
For each $t \in T$ with $v_{t} = v_{b}$ and $k > 1$, 
$\omega_{t, k}$ is a rational differential on $C_{\Delta} \otimes_{A_{\Delta}} B_{\Delta}$ 
which has only pole (of order $k$) at the point $p_{t}$ corresponding to $t$. 

{\rm (3)}
For each $t_{1}, t_{2} \in T$ such that $t_{1} \neq t_{2}$, 
$\omega_{t_{1}, t_{2}}$ is a rational differential on 
$C_{\Delta} \otimes_{A_{\Delta}} B_{\Delta}$ 
which has only (simple) poles at the points $p_{t_{1}}$ (resp. $p_{t_{2}}$) 
corresponding to $t_{1}$ (resp $t_{2}$) with residue $1$ (resp. $-1$). 

{\rm (4)} 
Take $x_{h}$ $(h \in \pm E \cup T)$ and $y_{e}$ $(e \in E)$ 
be complex numbers as in 2.2 (P4). 
Then $\omega_{i}, \omega_{t, k}, \omega_{t_{1}, t_{2}}$ are abelian differentials 
on the Riemann surface $R_{\Gamma}$, 
where $\Gamma$ is the Schottky group obtained from $\Gamma_{\Delta}$. 

\end{it} 
\vspace{2ex}

{\it Proof.} 
By Proposition 4.1, 
$\omega_{i}$ are differentials on ${\mathcal P}_{T_{\Delta}}$, 
and for any $\delta \in \Gamma_{\Delta}$, 

\begin{eqnarray*}
\omega_{i}(\delta(z)) & = & 
\sum_{\gamma \in \Gamma_{\Delta} / \left\langle \gamma_{i} \right\rangle} 
\left( \frac{\gamma(\alpha_{i}) - \gamma(\alpha_{-i})}
{(\delta(z) - \gamma(\alpha_{i})) (\delta(z) - \gamma(\alpha_{-i}))} \right) 
d \delta(z) 
\\ 
& = & 
\sum_{\gamma \in \Gamma_{\Delta} / \left\langle \gamma_{i} \right\rangle} 
\left( \frac{(\delta^{-1} \gamma)(\alpha_{i}) - (\delta^{-1} \gamma)(\alpha_{-i})}
{(z - (\delta^{-1} \gamma)(\alpha_{i})) (z - (\delta^{-1} \gamma)(\alpha_{-i}))} \right) dz 
\\ 
& = & 
\omega_{i}(z). 
\end{eqnarray*}
Therefore, by the construction of $C_{\Delta}$ reviewed in 2.2, 
$\omega_{i}$ give rise to differentials on $C_{\Delta}$ 
which are regular outside $\bigcup_{e \in E} \{ y_{e} = 0 \}$, 
and hence the assertions (1) follows. 
One can prove (2)--(3) similarly, 
and we prove (4). 
As is stated in 2.1, 
$R_{\Gamma}$ is given by the quotient space $\Omega_{\Gamma}/\Gamma$. 
Under the assumption on complex numbers $x_{h}$ and $y_{e}$, 
it is shown in \cite{S} that $\sum_{\gamma \in \Gamma} |\gamma'(z)|$ 
is uniformly convergent on any compact subset in 
$\Omega_{\Gamma} - \cup_{\gamma \in \Gamma} \gamma(\infty)$, 
and hence the assertion holds for $\omega_{t, k}$. 
If $a \in \Omega_{\Gamma} - \cup_{\gamma \in \Gamma} \gamma(\infty)$, 
then $\lim_{n \rightarrow \infty} \gamma_{i}^{\pm n}(a) = \alpha_{\pm i}$, 
and hence 
\begin{eqnarray*}
\lefteqn{
d \left( \int_{a}^{\gamma_{i}(a)} \sum_{\gamma \in \Gamma} 
\frac{d \gamma(\zeta)}{\gamma(\zeta) - z} \right)
} 
\\ 
& = & 
\sum_{\gamma \in \Gamma} 
\left( \frac{1}{z - (\gamma \gamma_{i})(a)} - \frac{1}{z - \gamma(a)} \right) dz 
\\ 
& = & 
\sum_{\gamma \in \Gamma / \langle \gamma_{i} \rangle} 
\lim_{n \rightarrow \infty} 
\left( \frac{1}{z - (\gamma \gamma_{i}^{n})(a)} - 
\frac{1}{z - (\gamma \gamma_{i}^{-n})(a)} \right) dz 
\\ 
& = & 
\omega_{i}(z). 
\end{eqnarray*}
Therefore, 
$\omega_{i}$ is absolutely and uniformly convergent 
on any compact subset in $\Omega_{\Gamma}$, 
and hence is an abelian differential on $\Omega_{\Gamma}/\Gamma$. 
\ $\square$ 
\vspace{2ex}

{\bf 4.2. Stability of differentials.} \ 
For a vertex $v \in V$, 
denote by $C_{v}$ the corresponding irreducible component of 
$C_{\Delta} \otimes_{A_{\Delta}} R_{\Delta}$. 
Then $P_{v} = {\mathbb P}^{1}_{R_{\Delta}}$ is the normalization of $C_{v}$. 
\vspace{2ex}

{\sc Theorem 4.3.} 
\begin{it} 

{\rm (1)} 
For each $i = 1,..., g$, 
let $\phi_{h_{i}(1)} \cdots \phi_{h_{i}(l_{i})}$ $(h_{i}(j) \in \pm E)$ be the unique reduced product 
such that $v_{-h_{i}(j)} = v_{h_{i}(j+1)}$ and $h_{i}(1) \neq - h_{i}(l_{i})$ which is conjugate to $\gamma_{i}$. 
Then for each $v \in V$, 
the pullback $\left( \omega_{i}|_{C_{v}} \right)^{*}$ of $\omega_{i}|_{C_{v}}$ to $P_{v}$ 
is given by 
$$
\left( \sum_{v_{h_{i}(j)} = v} \frac{1}{z - x_{h_{i}(j)}} - 
\sum_{v_{-h_{i}(k)} = v} \frac{1}{z - x_{-h_{i}(k)}} \right) dz. 
$$

{\rm (2)}
For each $v \in V$, 
$\left( \omega_{t, k}|_{C_{v}} \right)^{*}$ is given by 
$\displaystyle \frac{dz}{(z - x_{t})^{k}}$ if $v = v_{t}$, 
and is $0$ otherwise. 

{\rm (3)} 
Denote by $\rho_{t_{j}} = h_{j}(1) \cdots h_{j}(l_{j})$ 
the unique path from $v_{t_{j}}$ $(t_{j} \in T)$ to $v_{b}$ in ${\mathcal T}_{\Delta}$. 
Then for each $v \in V$, 
$\left( \omega_{t_{1}, t_{2}}|_{C_{v}} \right)^{*}$ is given by 
$$
\left( \sum_{v_{h} = v} \frac{1}{z - x_{h}} - 
\sum_{v_{-k} = v} \frac{1}{z - x_{-k}} \right) dz, 
$$
where $h, k$ runs through 
$\left\{ t_{1}, h_{1}(1),..., h_{1}(l_{1}), -t_{2}, -h_{2}(1),..., -h_{2}(l_{2}) \right\}$. 
\end{it} 
\vspace{2ex}

{\it Proof.} 
For the proof of (1), 
we may assume that 
$\gamma_{i} = \phi_{h_{i}(1)} \cdots \phi_{h_{i}(l_{i})}$. 
Let $\gamma$ be an element of $\Gamma_{\Delta}$. 
Then by Proposition 4.1 (2), putting $y_{e} = 0$ $(e \in E)$, 
$$
\frac{1}{z - \gamma(\alpha_{i})} - \frac{1}{z - \gamma(\alpha_{-i})} = 
\frac{\gamma(\alpha_{i}) - \gamma(\alpha_{-i})}
{(z - \gamma(\alpha_{i}))(z - \gamma(\alpha_{-i}))}  
$$
becomes 
$$
\frac{1}{z - x_{h_{i}(j)}} - \frac{1}{z - x_{-h_{i}(j-1)}}
$$ 
if $j \in \{ 1,..., l_{i} \}$ $\left( \mbox{$h_{i}(j-1) := h_{i}(l_{i})$ when $j = 1$} \right)$, 
$v(h_{i}(j)) = v$, 
$$
\phi_{h_{i}(j)} \phi_{h_{i}(j+1)} \cdots \phi_{h_{i}(l_{i})} \in \gamma \langle \gamma_{i} \rangle, 
$$
and becomes $0$ otherwise. 
Therefore, the assertion follows from the definition of $\omega_{i}$. 

The assertion (2) follows from Proposition 4.1 (1) and (3), 
and the assertion (3) can be shown in the same way as above. 
\ $\square$ 
\vspace{2ex}

{\it Remark} 4.4. 
By regarding $C_{\Delta}$ as a family of marked Riemann surfaces given in 2.2 (P4), 
for each $h \in \pm E$, 
let $c_{h}$ be its cycle corresponding to $-h$ which is oriented by the right-hand rule. 
Then under the notation in Theorem 4.3 (1) and (2), 
$$
\frac{1}{2 \pi \sqrt{-1}} \int_{c_{h}} \omega_{i} 
= \left\{ \begin{array}{ll} 
1   & \mbox{(if $h \in \{ h_{i}(j) \ | \ j = 1,..., l_{i} \}$),} 
\\ 
-1  & \mbox{(if $h \in \{ -h_{i}(j) \ | \ j = 1,..., l_{i} \}$),} 
\\ 
0   & \mbox{(otherwise)} 
\end{array} \right. 
$$ 
and 
$$
\frac{1}{2 \pi \sqrt{-1}} \int_{c_{h}} \omega_{t, k} 
= 0 
$$
respectively. 
Therefore, for $\omega_{i}$, $\omega_{t, k}$, 
one can describe their analytic characterizations and asymptotic behaviors 
under $y_{e} \rightarrow 0$ for some elements $e \in E$. 
Since marked Riemann surfaces obtained as in 2.2 (P4) make a nonempty open subset 
in the moduli space of marked Riemann surfaces, 
by the theorem of identity, 
similar asymptotic behaviors also holds 
for families of general marked Riemann surfaces (cf. \cite{I7}). 
This modification gives a more explicit formula than \cite[Corollary 4.6]{HN}. 

When $C_{\Delta}$ is regarded as a family of marked Mumford curves over 
a $p$-adic field given in 2.2 (P2), 
$\displaystyle \frac{1}{2 \pi \sqrt{-1}} \int_{c_{h}}$ can be replaced by 
the Schneider integral (cf. \cite[0.3]{dS2}). 
\vspace{2ex}

A  (regular or rational) global section of the dualizing sheaf on a stable curve 
is called a {\it stable differential} (cf. \cite{DM}). 
\vspace{2ex}

{\sc Theorem 4.5.} 
\begin{it} 

{\rm (1)} 
For each $i = 1,..., g$, 
$\omega_{i}$ is a regular stable differential on $C_{\Delta}$. 
Furthermore, $\{ \omega_{i} \}_{1 \leq i \leq g}$ gives a basis of 
$H^{0} \left( C_{\Delta}, \omega_{C_{\Delta}/A_{\Delta}} \right)$. 

{\rm (2)} 
For each $t \in T$ with $v_{t} = v_{b}$ and $k > 1$, 
$\omega_{t, k}$ is a rational stable differential on $C_{\Delta}$ 
which has only pole (of order $k$) at the point $p_{t}$ corresponding to $t$. 

{\rm (3)}
For each $t_{1}, t_{2} \in T$ with $t_{1} \neq t_{2}$, 
$\omega_{t_{1}, t_{2}}$ is a rational stable differential on $C_{\Delta}$ 
which has only (simple) poles at the points $p_{t_{1}}$ (resp. $p_{t_{2}}$) 
corresponding to $t_{1}$ (resp $t_{2}$) with residue $1$ (resp. $-1$). 

\end{it} 
\vspace{2ex}

{\it Proof.} 
We only show that the latter assertion in (1) 
since the remains follow from Theorems 4.2 and 4.3. 
Let $\omega$ be an element of $H^{0} \left( C_{\Delta}, \omega_{C_{\Delta}/A_{\Delta}} \right)$. 
Then by Remark 4.4, 
its residue ${\rm Res}_{h}(\omega) \in A_{\Delta}$ at $x_{h}$ $(h \in \pm E)$ satisfies 
$$
{\rm Res}_{h}(\omega) = \frac{1}{2 \pi \sqrt{-1}} \int_{c_{h}} \omega, 
$$ 
and ${\rm Res}_{h}(\omega) : \pm E \rightarrow A_{\Delta}$ gives 
a harmonic $1$-cochain in $\Delta$, namely 
$$
{\rm Res}_{h}(\omega) = - {\rm Res}_{-h}(\omega) \ (h \in \pm E), \ \ 
\sum_{v_{h} = v} {\rm Res}_{h}(\omega) = 0 \ (v \in V). 
$$
By the map sending harmonic $1$-cochains in $\Delta$ to 
the associated singular cohomology classes, 
the space of $A_{\Delta}$-valued harmonic $1$-cochain in $\Delta$ 
is canonically isomorphic to $H^{1}(\Delta, A_{\Delta})$. 
Since $\{ \gamma_{i} \}_{1 \leq i \leq g}$ gives a basis of 
$H_{1}(\Delta, {\mathbb Z}) \cong 
\Gamma_{\Delta}/\left[ \Gamma_{\Delta}, \Gamma_{\Delta} \right]$, 
there exist uniquely $r_{i} \in A_{\Delta}$ $(i = 1,..., g)$ such that 
$$
{\rm Res}_{h} \left( \omega - \sum_{i = 1}^{g} r_{i} \omega_{i} \right) = 0 \ \ (h \in \pm E) 
$$ 
which implies that $\omega = \sum_{i = 1}^{g} r_{i} \omega_{i}$. 
\ $\square$ 
\vspace{2ex}

{\it Remark} 4.6. 
By the proof of Theorem 4.5 (1), for each $i = 1,..., g$, 
the differential $\omega_{i}$ of the first kind is uniquely determined by the element 
$[ \gamma_{i} ]$ of $H_{1}({\Delta}, {\mathbb Z})$ derived from $\gamma_{i}$.  
\vspace{4ex}

{\bf 5. Universal differentials and periods.} 
\vspace{2ex}

{\bf 5.1. Universal differentials.} \ 
Denote by $\Delta_{0} = (V_{0}, E_{0}, T_{0})$ the stable graph of $(g, n)$-type 
consisting of one vertex and $g$ loops, 
and fix generators $\rho_{1},..., \rho_{g}$ of $\pi_{1}(\Delta_{0})$. 
For each stable graph $\Delta = (V, E, T)$ of $(g, n)$-type, 
$\Delta_{0}$ is obtained from $\Delta$ by contracting some nonloop edges in $E$. 
Then there exists a unique ${\mathbb Z}$-basis $\{ b_{i} \}_{1 \leq i \leq g}$ 
of $H_{1}(\Delta, {\mathbb Z})$ corresponding to 
$$
\left\{ [\rho_{i}] \right\}_{1 \leq i \leq g} \in H_{1}(\Delta_{0}, {\mathbb Z}), 
$$
and $T$ is identified with $T_{0}$.  
Therefore, by Theorem 4.5 and Remark 4.6, 
there exist associated stable differentials 
$$
\omega_{i} \ (i = 1,..., g), \ \ \omega_{t, k} \ (t \in T_{0}, k > 1), \ \ 
\omega_{t_{1}, t_{2}} \ (t_{i} \in T_{0}, t_{1} \neq t_{2}) 
$$
on $C_{\Delta}$ of the first, second, third kind respectively. 
Furthermore, the above $b_{i} \in H_{1}(\Delta, {\mathbb Z})$ $(i = 1,..., g)$ 
give rise to homology cycles on families of Riemann surfaces obtained from $C_{\Delta}$ 
as in 2.2 (P4), 
and on the universal family of stable complex curves over 
the extended Schottky space $\overline{S}_{g,n}$ by the analytic extension. 
We also denote by $b_{i}$ these homology cycles. 
\vspace{2ex}

{\sc Theorem 5.1.} 
\begin{it} 

{\rm (1)} 
The differentials $\omega_{i}$, $\omega_{t, k}$, $\omega_{t_{1}, t_{2}}$ on $C_{\Delta}$ 
are glued to stable differentials on ${\mathcal C}_{g}/{\mathcal A}_{g,n}$ 
which we call the universal differentials of the first kind, second kind, third kind, 
and denote by 
$\overline{\omega}_{i}$, $\overline{\omega}_{t, k}$, $\overline{\omega}_{t_{1}, t_{2}}$ respectively. 

{\rm (2)} 
The universal differentials $\overline{\omega}_{i}$ $(i = 1,..., g)$ on 
${\mathcal C}_{g}/{\mathcal A}_{g,n}$ of the first kind make a basis of the sheaf 
${\mathcal H}^{0} \left( \omega_{{\mathcal C}_{g}/{\mathcal A}_{g, n}} \right)$ 
which consists of sections of the relative stable differentials 
on ${\mathcal C}_{g}/{\mathcal A}_{g,n}$. 
Furthermore, $\overline{\omega}_{i}$ are analytically extended to 
stable differentials on the universal family of stable complex curves over $\overline{S}_{g,n}$ 
which we denote by the same symbols. 
\end{it}
\vspace{2ex} 

{\it Proof.} 
First, we prove (1). 
As is stated in Remark 3.2, 
the isomorphism $C_{\Delta} \cong C_{\Delta_{0}}$ over $B_{\Delta}$ 
considered in Theorem 3.1 corresponds uniquely to an isomorphism 
$\Gamma_{\Delta} \cong \Gamma_{\Delta_{0}}$. 
Therefore, under this isomorphism $C_{\Delta} \cong C_{\Delta_{0}}$, 
the differentials $\omega_{i}$, $\omega_{t, k}$, $\omega_{t_{1}, t_{2}}$ on $C_{\Delta}$ 
are mapped to those on $C_{\Delta_{0}}$, 
and hence can be glued to differentials on ${\mathcal C}_{g}/{\mathcal A}_{g,n}$.   

Second, we prove the latter assertion of (2) 
since the former one follows from Theorem 4.5 (1). 
For a stable $n$-marked complex curve $C$ of genus $g$ with Schottky structure, 
as is stated in \cite{GH}, 
one can take a cut system $\{ a_{1},..., a_{g} \}$ consisting of disjoint oriented simple loops 
in $C$ such that the intersection numbers $a_{i} \cdot b_{j}$ 
are the Kronecker delta $\delta_{ij}$. 
Then there exist uniquely regular stable differentials $\omega_{C, i}$ $(i = 1,..., g)$ on $C$ 
such that 
$$
\displaystyle \int_{a_{j}} \omega_{C, i} = 2 \pi \sqrt{-1} \delta_{ij}. 
$$
Therefore, by Theorem 4.3 (1), 
moving $C$ on $\overline{S}_{g, n}$ 
$\omega_{C, i}$ form the analytic extension of $\overline{\omega}_{i}$.  
\ $\square$
\vspace{2ex}

{\bf 5.2. Universal periods and Jacobian.} \ 
Denote by ${\mathcal B}_{g, n}$ the maximal open subscheme of  ${\mathcal A}_{g, n}$ 
over which the associated marked curves are smooth, 
and by $S_{g, n}$ the (ordinary) Schottky space 
which is the open subspace of $\overline{S}_{g, n}$ 
classifying $n$-marked Riemann surfaces of genus $g$ with Schottky structure. 
\vspace{2ex}

{\sc Theorem 5.2} 
\begin{it} 
There exist ${\mathcal P}_{ij} \in {\mathcal O}^{\times}_{{\mathcal B}_{g, n}}$ 
$(1 \leq i, j \leq g)$ which give the multiplicative periods of Mumford curves 
over $p$-adic fields obtained from ${\mathcal C}_{g}/{\mathcal B}_{g, n}$ as in 2.2 (P2). 
Furthermore, ${\mathcal P}_{ij}$ are analytically continued to regular functions 
on $S_{g, n}$ which are $\displaystyle \exp \left( \int_{b_{i}} \overline{\omega}_{j} \right)$, 
where $\overline{\omega}_{j}$ are given in Theorem 5.1 (2). 
We call ${\mathcal P}_{ij}$ the universal periods. 
\end{it}
\vspace{2ex}

{\it Proof.} 
For a stable graph $\Delta$ of $(g, n)$-type, 
take generators $\gamma_{1},..., \gamma_{g}$ of $\Gamma_{\Delta}$ such that 
$[\gamma_{i}] = b_{i}$ $(i = 1,..., g)$, 
and for each $\gamma_{i}$, 
denote by $\alpha_{i}$ (resp. $\alpha_{-i}$) its attractive (resp. repulsive) fixed points 
and by $\beta_{i}$ its multiplier. 
Then by \cite[Theorem 3.13]{I1}, 
the multiplicative periods $P_{ij}$ $(1 \leq i, j \leq g)$ of $C_{\Delta}$ are 
elements of $B_{\Delta}^{\times}$ defined as $P_{ij} = \prod_{\gamma} \psi_{ij}(\gamma)$, 
where $\gamma$ runs through all representatives of 
$\langle \gamma_{i} \rangle \backslash \Gamma_{\Delta} / \langle \gamma_{j} \rangle$ 
and 
$$
\psi_{ij}(\gamma) = \left\{ \begin{array}{ll} 
\beta_{i} & (i = j, \ \gamma \in \langle \gamma_{i} \rangle), 
\\
{\displaystyle \frac{(\alpha_{i} - \gamma(\alpha_{j})) (\alpha_{-i} - \gamma(\alpha_{-j}))}
{(\alpha_{i} - \gamma(\alpha_{-j})) (\alpha_{-i} - \gamma(\alpha_{j}))}} 
& (\mbox{otherwise}). 
\end{array} \right. 
$$ 
By Remark 4.6, $P_{ij}$ depend only on $\{ b_{i} \}_{1 \leq i \leq g}$, 
and one can see that they give rise to the multiplicative periods of Mumford curves 
over $p$-adic fields by \cite[Theorem 2]{MD} 
and to $\displaystyle \exp \left( \int_{b_{i}} \omega_{j} \right)$ 
by the definition of $\omega_{j}$ (cf. \cite{S}). 
Therefore, as in the proof of Theorem 5.1, 
$P_{ij}$ are seen to be glued to regular functions ${\mathcal P}_{ij}$ on ${\mathcal B}_{g, n}$ 
satisfying the required properties. 
\ $\square$ 
\vspace{2ex}

{\sc Theorem 5.3.} 
\begin{it} 
There exists a polarized semi-abelian scheme ${\mathcal J}_{g}$ over ${\mathcal A}_{g, n}$ 
whose fibers at ${\rm Spec}(A_{\Delta})$ ($\Delta$: stable graphs of $(g, n)$-type) 
are the polarized semi-abelian scheme $J_{\Delta}$ given in \cite[Theorem 3.13]{I1} 
as $J_{\Pi}$. 
Consequently, ${\mathcal J}_{g}$ gives rise to the Jacobian varieties of 
Riemann surfaces and of Mumford curves obtained from 
${\mathcal C}_{g}/{\mathcal A}_{g, n}$ as in 2.2 (P4) and (P2) respectively. 
Furthermore, ${\mathcal J}_{g}$ is analytically extended to a family of 
$g$-dimensional semi-abelien varieties over the quotient space 
$\overline{S}_{g, n}/{\rm Out}(F_{g})$. 
We call ${\mathcal J}_{g}$ the universal Jacobian. 
\end{it}
\vspace{2ex}

{\it Proof.} 
First, we recall the construction of $J_{\Delta}$ for a stable graph $\Delta$ of $(g, n)$-type. 
Put $X = H_{1}(\Delta, {\mathbb Z})$, 
and let ${\mathbb G}_{m}^{g}$ denote the split torus of dimension $g$ over $A_{\Delta}$ 
with character group $X$. 
Let $L$ be the group of periods of $C_{\Delta}$ which is defined as the subgroup 
of $B_{\Delta}^{\times}$ generated by the above $(P_{ij})_{1 \leq i \leq g}$, namely 
$$
L = \left\{ \left. \left( \prod_{j=1}^{g} P_{ij}^{n_{j}} \right)_{1 \leq i \leq g} \ \right| 
n_{j} \in {\mathbb Z} \right\}. 
$$
Then the isomorphism $\phi: L \rightarrow X$ given by 
$$
\phi \left( \left( \prod_{j=1}^{g} P_{ij}^{n_{j}} \right)_{1 \leq i \leq g} \right) ((z_{i})_{i}) 
= \prod_{i=1}^{g} z_{i}^{n_{i}} \ \left( (z_{i})_{i} \in {\mathbb G}_{m}^{g} \right) 
$$ 
satisfies the definition of polarizations in \cite[Definition 1.2]{Mu2}, 
and hence by results of \cite{Mu2} and \cite[Chapter III]{FC}, 
there exists a semi-abelian scheme $J_{\Delta}$ over the regular ring $A_{\Delta}$ 
which is formally represented as the quotient of ${\mathbb G}_{m}^{g}$ by $L$. 
By this representation, for each integer $r \geq 1$, 
the subgroup of ${\mathbb G}_{m}^{g}$ of $r$th roots of $1$ gives a canonical subgroup 
$\mbox{\boldmath $\mu$}_{r}$ of $J_{\Delta}$. 

From the above construction of $J_{\Delta}$, 
the above isomorphism $C_{\Delta} \cong C_{\Delta_{0}}$ 
gives rise to a unique isomorphism $J_{\Delta} \cong J_{\Delta_{0}}$ 
which is the identity map on $\mbox{\boldmath $\mu$}_{r}$ for $r \geq 3$. 
Therefore, by Theorem 5.2, 
$J_{\Delta}$ ($\Delta$: stable graphs of $(g, n)$-type) are glued to 
a semi-abelian scheme ${\mathcal J}_{g}$ over ${\mathcal A}_{g,n}$ 
satisfying the required properties except the last assertion which we will prove. 
For a stable complex curve $C$ of genus $g$ with Schottky structure, 
take a stable curve ${\mathcal C}$ over a nonarchimedean valuation ring $R$ 
such that its generic fiber ${\mathcal C}_{\eta}$ is smooth and its special fiber gives $C$. 
Then the dual graph $\Delta_{C}$ of $C$ can be regarded as a subgraph of 
a stable graph of $(g, 0)$-type, 
and (after possible base change) the connected Neron model ${\rm Jac}^{\circ}({\mathcal C})$ 
of the Jacobian variety of ${\mathcal C}_{\eta}$ is a semi-abelian scheme over $R$ 
whose torus part has the character group $H_{1}(\Delta_{C}, {\mathbb Z})$. 
Moving $C$ over $\overline{S}_{g}/{\rm Out}(F_{g})$, 
the special fibers of ${\rm Jac}^{\circ}({\mathcal C})$ form a family of 
$g$-dimensional semi-abelien varieties which are derived from ${\mathcal J}_{g}$ 
if $C$ are close to degenerate curves by the functoriality of the Mumford construction 
of semi-abelian schemes \cite[Chapters III and VI]{FC}.    
\ $\square$
\vspace{2ex}

{\sc Theorem 5.4.} 
\begin{it} 
Assume that $n > 1$. 
Then for each $i = 1,..., g$, $t \in T_{0}$ and $k > 1$, 
$\displaystyle \int_{b_{i}} \overline{\omega}_{t, k}$ can be defined as 
a regular function on ${\mathcal A}_{g, n} \otimes {\mathbb Q}$ which gives 
the associated line integral for Riemann surfaces 
and the Coleman integral for Mumford curves over $p$-adic fields 
obtained from ${\mathcal C}_{g, n}/{\mathcal A}_{g, n}$ as in 2.2. 
\end{it}
\vspace{2ex}

{\it Proof.} 
For a stable graph $\Delta$ of $(g, n)$-type, 
take $\gamma_{i} \in \Gamma_{\Delta}$ $(i = 1,..., g)$ such that $[\gamma_{i}] = b_{i}$. 
By the assumption, there is an element $t' \in T_{0} - \{ t \}$, 
and hence one can take a path $\rho_{t'}$ in $\Delta$ from $v_{t'}$ to $v_{t}$ 
and an element $\phi_{t'}$ of $\Gamma_{\Delta}$ corresponding to $\rho_{t'}$. 
Then by Proposition 4.1, 
$$
(1 - k) \int_{b_{i}} \omega_{t, k} = \sum_{\gamma \in \Gamma_{\Delta}} 
\left( \frac{1}{(\gamma(\gamma_{i}(\phi_{t'}))) - 
x_{t})^{k-1}} - \frac{1}{(\gamma(\phi_{t'})) - x_{t})^{k-1}} \right)
$$ 
is defined as an element of $A_{\Delta}$, 
and $\displaystyle \int_{b_{i}} \omega_{t, k} \in A_{\Delta} \otimes {\mathbb Q}$ 
($\Delta$: stable graphs of $(g, n)$-type) are glued to a regular function 
on ${\mathcal A}_{g, n} \otimes {\mathbb Q}$ 
which we denote by $\displaystyle \int_{b_{i}} \overline{\omega}_{t, k}$. 
By results of de Shalit (cf. \cite[1.5]{dS1}, \cite[0.4]{dS2}), 
this function satisfies the required properties. 
\ $\square$ 
\vspace{2ex}

Denote by ${\mathcal H}_{\rm dR}^{1} \left( {\mathcal C}_{g}/{\mathcal B}_{g, n} \right)$ 
the sheaf of the first relative de Rham cohomology groups 
of ${\mathcal C}_{g}/{\mathcal B}_{g, n}$ which contains 
${\mathcal H}^{0} \left( \Omega_{{\mathcal C}_{g}/{\mathcal B}_{g, n}} \right)$ as its subsheaf. 
\vspace{2ex}

{\sc Theorem 5.5.} 
\begin{it} 
Assume that $n > 1$ and there exist differentials $\overline{\omega}_{t_{j}, k_{j}}$ 
$(j = 1,..., g)$ of the second kind on ${\mathcal A}_{g, n}$ such that 
$$
d \left( \overline{\omega}_{t_{j}, k_{j}} \right) := 
\det \left( \int_{b_{i}} \overline{\omega}_{t_{j}, k_{j}} \right)_{1 \leq i, j \leq g} 
$$
is a nonzero function on ${\mathcal A}_{g, n} \otimes {\mathbb Q}$. 
We denote by $U$ the maximal open subscheme of 
${\mathcal B}_{g, n} \otimes {\mathbb Q}$ such that 
$d \left( \overline{\omega}_{t_{j}, k_{j}} \right) \in {\mathcal O}^{\times}_{U}$. 
Then there exist sections $\overline{\eta}_{i}$ $(i = 1,..., g)$ of 
${\mathcal H}_{\rm dR}^{1} \left( {\mathcal C}_{g}/U \right)$ such that 
$\left\{ \overline{\omega}_{i}, \overline{\eta}_{i} \right\}_{1 \leq i \leq g}$ gives 
a basis of ${\mathcal H}_{\rm dR}^{1} \left( {\mathcal C}_{g}/U \right)$ 
and that the Gauss-Manin connection 
$$
\nabla : {\mathcal H}_{\rm dR}^{1} \left( {\mathcal C}_{g}/{\mathcal B}_{g, n} \right) 
\rightarrow  {\mathcal H}_{\rm dR}^{1} \left( {\mathcal C}_{g}/{\mathcal B}_{g, n} \right) 
\otimes \Omega_{{\mathcal B}_{g, n}} 
$$
satisfies 
$\displaystyle \nabla( \overline{\omega}_{i} ) = 
\sum_{j = 1}^{g} \overline{\eta}_{j} \otimes \left( d {\mathcal P}_{ij}/{\mathcal P}_{ij} \right)$, 
$\nabla( \overline{\eta}_{i} ) = 0$ on $U$.
\end{it}
\vspace{2ex}

{\it Proof.} 
By the definition of $U$, 
there exist differentials $\overline{\eta}'_{j}$ $(j = 1,..., g)$ given as 
${\mathcal O}_{U}$-linear sums of $\overline{\omega}_{t_{j}, k_{j}}$ such that 
$\displaystyle \int_{b_{i}} \overline{\eta}'_{j} = \delta_{ij}$. 
Then by Remark 4.4, 
$\displaystyle \int_{c_{h}} \overline{\eta}'_{j} = 0$ for any $h \in \pm E$, 
and hence $\overline{\eta}'_{j}$ give rise to sections $\overline{\eta}_{j}$ of 
${\mathcal H}_{\rm dR}^{1} \left( {\mathcal C}_{g}/U \right)$ which are dual to  $\overline{\omega}_{i}$ for the canonical symplectic form on 
$$
{\mathcal H}_{\rm dR}^{1} \left( {\mathcal C}_{g}/U \right) \times 
{\mathcal H}_{\rm dR}^{1} \left( {\mathcal C}_{g}/U \right). 
$$
Therefore, $\left\{ \overline{\omega}_{i}, \overline{\eta}_{i} \right\}_{1 \leq i \leq g}$ 
satisfies the required properties. 
\ $\square$  
\vspace{4ex}

\begin{flushleft} 
{\sc Department of Mathematics, Faculty of Science and Engineering, 
Saga University, Saga 840-8502, Japan} 
\\
{\it E-mail:}  ichikawn@cc.saga-u.ac.jp 
\end{flushleft} 
\vspace{2ex}

\hspace{4cm} \hrulefill \hspace{4cm}

\renewcommand{\refname}{\centerline{\scriptsize{REFERENCES}}}
\bibliographystyle{amsplain}

\begin{thebibliography}{99}


\bibitem{A} 
V. Alexeev, 
Compactified Jacobians and Torelli map, 
{\it Publ. RIMS, Kyoto Univ.} {\bf 40} (2004), 1241--1265. 

\bibitem{BM} 
M. Bainbridge and M. M\"{o}ller, 
The Deligne-Mumford compactification of the real multiplication locus and 
Teichm\"{u}ller curves in genus $3$, 
{\it Acta Math.} {\bf 208} (2012), 1--92. 

\bibitem{BaK1}  
B. Bakalov and A. Kirillov, 
On the Lego-Teichm\"{u}ller game, 
{\it Transform. Groups} {\bf  5} (2000), 207--244.  

\bibitem{BaK2}  
B. Bakalov and A. Kirillov,   
Lectures on Tensor categories and modular functors, 
{\it University Lecture Series, vol. 21,} Amer. Math. Soc, 2001, 
Available at http://www.math.stonybrook.edu/~kirillov/tensor/tensor.html

\bibitem{C} 
L. Caporaso, 
A compactification of the universal Picard variety over the moduli space of stable curves, 
{\it J. Amer. Math. Soc.} {\bf 7} (1994), 589--660. 

\bibitem{ChPS} 
K. Christ, S. Payne and J. Shen, 
Compactified Jacobians as Mumford models, 
arXiv:1912.03653. 

\bibitem{DM} 
P. Deligne and D. Mumford, 
The irreducibility of the space of curves of given genus, 
{\it Inst. Hautes \'{E}tudes Sci. Publ. Math.} {\bf 36} (1969), 75--109. 

\bibitem{dS1} 
E. de Shalit, 
Differentials of the second kind on Mumford curves, 
{\it Israel J. of Math.} {\bf 71} (1990), 1--16. 

\bibitem{dS2} 
E. de Shalit, 
Coleman integration versus Schneider integration on semistable curves, 
{\it Doc. Math. Extra Volume Coates} (2006), 325--334. 

\bibitem{E}
E. Esteves, 
Compactifying the relative Jacobian over families of reduced curves, 
{\it Trans. Amer. Math. Soc.} {\bf 353} (2001), 3045--3095. 

\bibitem{FC} 
G. Faltings and C. L. Chai, 
Degeneration of abelian varieties, 
{\it Ergeb. Math. Grenzgeb., vol. 22,} Springer-Verlag, Berlin 1990. 

\bibitem{GH} 
L. Gerritzen and F. Herrlich, 
The extended Schottky space, 
{\it J. reine angew. Math.} {\bf 389} (1988), 190--208. 

\bibitem{Gr} 
A. Grothendieck, 
Esquisse d'un programme. Mimeographed Note (1984), 
{\it Geometric Galois action I, London Math. Soc. Lect. Note Ser., vol. 242,} 
London Math. Soc., 1997, pp. 5--48. 

\bibitem{HN} 
X. Hu and C. Norton, 
General variational formulas for abelian differentials, 
{\it Int. Math. Res. Not.} {\bf 2020} (2020), 3540--3581. 

\bibitem{I1} 
T. Ichikawa, 
Generalized Tate curve and integral Teichm\"{u}ller modular forms, 
{\it Amer. J. Math.} {\bf 122} (2000), 1139--1174. 

\bibitem{I2} 
T. Ichikawa, 
Teichm\"{u}ller groupoids and Galois action, 
{\it J. reine angew. Math.} {\bf 559} (2003), 95--114. 

\bibitem{I3} 
T. Ichikawa, 
Teichm\"{u}ller groupoids, and monodromy in conformal field theory, 
{\it Commun. Math. Phys.} {\bf 246} (2004), 1--18. 

\bibitem{I4}
T. Ichikawa, 
Klein's formulas and arithmetic of Teichm\"{u}uller modular forms, 
{\it Proc. Amer. Math. Soc.} {\bf 146} (2018), 5105--5112. 

\bibitem{I5}
T. Ichikawa, 
Chern-Simons invariant and Deligne-Riemann-Roch isomorphism, 
{\it Trans. Amer. Math. Soc.} {\bf 374} (2021), 2987--3005. 

\bibitem{I6} 
T. Ichikawa, 
An explicit formula of the normalized Mumford form, 
{\it Lett. Math. Phys.} (2021) 111:2, 
https://doi.org/10.1007/s11005-020-01339-0. 

\bibitem{I7} 
T. Ichikawa, 
Stable degeneration of abelian differentials and of quasi-periodic solutions 
of the KP hierarchy, 
Submitted. 

\bibitem{I8}
T. Ichikawa, 
The universal Mumford curve and its periods in arithmetic formal geometry, 
arXiv:2010.11517v2. 

\bibitem{IhN} 
Y. Ihara and H. Nakamura, 
On deformation of maximally degenerate stable marked curves and Oda's problem, 
{\it J. reine angew. Math.} {\bf 487} (1997), 125--151. 

\bibitem{K} 
F. F. Knudsen, 
The projectivity of the moduli space of stable curves II, III, 
{\it Math. Scand.} {\bf 52} (1983), 161--199, 200--212. 

\bibitem{KM} 
F. F. Knudsen and D. Mumford, 
The projectivity of the moduli space of stable curves I, 
{\it Math. Scand.} {\bf 39} (1976), 19--55. 

\bibitem{MD} 
Yu. Manin and V. Drinfeld, 
Periods of $p$-adic Schottky groups, 
{\it J. Reine Angew. Math.} {\bf 262/263} (1972), 239--247.  

\bibitem{MoS} 
G. Moore and N. Seiberg, 
Classical and quantum conformal field theory, 
{\it Commun. Math. Phys.} {\bf 123} (1989), 177--254. 

\bibitem{Mu1} 
D. Mumford, 
An analytic construction of degenerating curves over complete local rings, 
{\it Compos. Math.} {\bf 24} (1972), 129--174. 

\bibitem{Mu2} 
D. Mumford, 
An analytic construction of degenerating abelian varieties over complete rings, 
{\it Compos. Math.} {\bf 24} (1972), 239--272. 

\bibitem{NS} 
H. Nakamura and L. Schneps, 
On a subgroup of Grothendieck-Teichm\"{u}ller group acting on the tower of 
profinite Teichm\"{u}ller modular groups, 
{\it Invent. Math.} {\bf 141} (2000), 503--560. 

\bibitem{OS} 
T. Oda and C. S. Seshadri, 
Compactifications of the generalized jacobian variety, 
{\it Trans. Amer. Math. Soc.} {\bf 253} (1979), 1--90. 

\bibitem{P1}
J. Poineau, 
La droite de Berkovich sur ${\mathbb Z}$, 
{\it Ast\'{e}risque} {\bf 334} (2010), xii+284. 

\bibitem{P2} 
J. Poineau, 
Espaces de Berkovich sur ${\mathbb Z}$: \'{e}tude locale. 
{\it Invent. Math.} {\bf 194} (2013), 535--590. 

\bibitem{PT} 
J. Poineau and D. Turchetti, 
Schottky spaces and universal Mumford curves over ${\mathbb Z}$, 
https://poineau.users.lmno.cnrs.fr/Textes/MumfordZ.pdf. 

\bibitem{R} 
M. Raynaud, 
G\'{e}om\'{e}trie analytique rigide d'apr\`{e}s Tate, Kiehl... 
{\it M\'{e}moires de la S. M. F.} tome 39--40 (1974), 319--327. 

\bibitem{S} 
F. Schottky, 
\"{U}ber eine specielle Function, welche bei einer bestimmten linearen Transformation ihres Arguments unver\"{a}ndert bleibt, 
{\it J. reine angew. Math.} {\bf 101} (1887), 227--272. 

\bibitem{Si} 
C. Simpson, 
Moduli of representations of the fundamental group of a smooth projective variety. I, 
{\it Inst. Hautes \'{E}tudes Sci. Publ. Math.} {\bf 79} (1994), 47--129. 

\bibitem{U} 
M. Ulirsch, 
A non-archimedean analogue of Teichm\"{u}ller space and its tropicalization, 
arXiv:2004.07508.  

\end{thebibliography}

\end{document}